\documentclass[10pt]{article}

\usepackage{amssymb}
\usepackage{latexsym}
\usepackage{amsmath}
\usepackage{amscd}
\usepackage{comb}
\usepackage{amsbsy}

\numberwithin{equation}{section}

\newtheorem{theorem}{Theorem}

\newtheorem{lemma}{Lemma}

\begin{document}
\title{Minimax estimation of the Wigner function in quantum homodyne tomography with ideal detectors}
\author{M\u ad\u alin Gu\c{t}\u{a} \footnote{Mathematical Institute, University of Utrecht, Budapestlaan 6, 3584 CD Utrecht, The Netherlands, guta@math.uu.nl}
\hspace{3mm} and  \hspace{3mm}
Luis  Artiles\footnote{Eurandom, University of Eindhoven, P.O. Box 513, 5600 MB Eindhoven, The Netherlands,
artiles@eurandom.tue.nl, http://euridice.tue.nl/~lartiles/}
}
\date{}
\maketitle


\begin{abstract}

We estimate the quantum state of a light beam from results of quantum
homodyne measurements performed on identically prepared pulses. 
The state is represented through the Wigner function, a ``quasi-probability density'' on 
$\mathbb{R}^{2}$ which may take negative values and must respect intrinsic
positivity constraints imposed by quantum physics. The data consists of $n$ i.i.d. observations from a probability density equal to the Radon transform of the Wigner function. We construct an estimator for the Wigner function, and prove that it is
minimax efficient for the pointwise risk over a class of infinitely 
differentiable functions. A similar result was previously derived by Cavalier in the 
context of positron emission tomography. Our work extends this result to 
the space of smooth Wigner functions, which is  the relevant parameter space 
for quantum homodyne tomography.


\end{abstract}

\section{Introduction}
\label{sec.introduction}

The phenomena occurring at the interface between the microscopic and the macroscopic 
worlds have an intrinsic probabilistic nature. When measuring properties of atoms and laser pulses we obtain a random result whose probability distribution is determined by the state, or preparation of the quantum system.  For example, if we count the number of photons coming from a laser source we observe a Poisson distributed random variable with mean equal to the intensity of the laser. 
The statistical inverse problem of inferring  the state from results of measurements on many identically prepared systems, is called quantum state estimation. Recently it has become possible to apply such a method to the reconstruction of the quantum state of a light beam. The measurement technique is called Quantum Homodyne Tomography \cite {Vogel&Risken} and is used  to confirm the creation of new and exotic quantum states of light such as squeezed states \cite{Breitenbach&Schiller&Mlynek}, single-photon-added coherent states  \cite{Zavatta} and Schr\" odinger cat states \cite{Ourjoumtsev}. As experiments become more and more complicated, the costs -- in terms of money and time -- of running a measurement rise, and one needs to apply more sophisticated statistical techniques to reconstruct the state from a limited number of data. 
This paper makes a step in this direction by providing minimax convergence rates for a class of physical states.

The object to be estimated is a real function of two variables called the 
Wigner function \cite{Wigner}, which can be seen as a joint density of the electric and magnetic fields of the laser beam. However, since in quantum mechanics we cannot measure both electric and magnetic fields simultaneously, this function is in general not a probability density but has many features in common with the latter, for example the marginals along any direction are bona-fide probability densities. In quantum optics 
the Wigner function is a preferred representation of the quantum state \cite{Glauber} because many interesting quantum patterns such as squeezing or oscillation between negative and positive values, can be easily identified from the shape of the function. The estimation methods used by the physicists involve a number of ad-hoc approximations, 
binnings and truncations, making it difficult to verify the reliability of the procedure. Moreover, the quantum features in which the experimenter is interested may be 
washed out in the resulting estimator.

From the statistical viewpoint, we deal with an ill posed inverse problem which is closely related to the problem of estimating a bivariate probability density in  Positron Emission Tomography. In both cases the parameter of interest is a density and the data consists of independent identically distributed samples from the Radon Transform \cite{Deans} of that density, with uniformly distributed angles. However as we mentioned above, the Wigner function is not necessarily positive but must satisfy other positivity conditions dictated by the laws of quantum mechanics. 
This means that some of the existing statistical results are not 
directly applicable in this case, and the problem of estimating the Wigner function should be studied separately in order to identify the specific ``quantum features'' that could be exploited 
in designing new estimators. The problem of state estimation in Quantum Homodyne Tomography was investigated in the non-parametric setting in \cite{Artiles&Gill&Guta}. The paper provides sufficient conditions for consistency of various estimators for two different representations of the state, namely the density matrix and the Wigner function. The problem of estimating the Wigner function was further investigated in 
\cite{Butucea&Guta&Artiles} in a set-up which takes into account the detection losses occurring in the measurement, leading to an additional additive Gaussian noise. 
Although both \cite{Butucea&Guta&Artiles} and the present paper consider a family of very smooth Wigner functions, the estimation techniques are very different. 
In \cite{Butucea&Guta&Artiles} the bias dominates the variance due to the presence of the Gaussian noise, and for deriving the lower bound it is enough to consider a ``worst family'' consisting of just two states. In this paper, the variance is dominating and for the lower bound we need to consider a one parameter family of states. The techniques that we use here are very similar to the one of Cavalier \cite{Cavalier} who considers the same tomography  problem in the context of smooth probability densities rather than Wigner functions. For the problem of  Positron Emission Tomography and the related regression problem of X-ray Tomography we refer to 
\cite{Johnstone.Silverman.1,Korostelev&Tsybakov,Cavalier02} and the references therein. For an introduction to Quantum Statistical Inference we refer to 
\cite{Barndorff-Nielsen&Gill&Jupp,Artiles&Gill&Guta}. 

In Section \ref{sec.2} we give a short introduction to quantum mechanics, present some  properties of Wigner functions and the relation to the density matrix which will play an important role in the proof of the lower bound. The statistical set-up and the main results are presented in Section \ref{sec.stat.setup}. Following \cite{Cavalier}, we identify a class of very smooth Wigner functions
\beq
\mathcal{W} (\beta, L) = \left\{ W  : 
\frac{1}{4\pi^{2}}\int_{\R^{2}} |\widetilde W (w)|^{2} \exp(2\beta |w|) \, dw \leq L \right\},
\eeq
where $\beta$ and $L$ are positive constants. From the physical point of view, all the states which have been produced in the lab up to date belong to such a class, 
and a more detailed argument as to why this assumption is realistic, can be found in \cite{Butucea&Guta&Artiles}.

We consider a family of estimators depending on a bandwidth which is chosen according to the parameters of the class. The upper bound for the pointwise risk is proven in Theorem \ref{theorem.upperbound} and has the same almost parametric expression as the bound derived in \cite{Cavalier}. For the lower bound we consider a family of 
Wigner functions which is different from the worst parametric family of probability densities of \cite{Cavalier}. The latter cannot be used in our situation since it does 
not correspond to Wigner functions of physical states.  Thus the main novelty of the paper is the derivation of the lower bound in the physical context of Quantum Homodyne Tomography rather than that of Positron Emission Tomography.

The technical Lemmas are grouped in Section \ref{sec.lemmas}.

\section{Quantum Homodyne Tomography and the Wigner function}\label{sec.2}

In this section we briefly present some basic notions of quantum mechanics, the mathematical set-up of Quantum Homodyne Tomography, and some properties of 
Wigner functions. More information on the general set-up of quantum statistical 
inference can be found in  the review paper \cite{Barndorff-Nielsen&Gill&Jupp} and the textbooks \cite{Helstrom} and \cite{Holevo}. Quantum Homodyne Tomography 
is discussed in detail in \cite{Artiles&Gill&Guta,Butucea&Guta&Artiles}.

The state of a quantum system encodes all necessary information for 
computing the probability distribution of results of any given measurement. Mathematically, the state is described by a {\it density matrix}, which is a compact
operator $\rho$ on a complex Hilbert space $\mathcal{H}$ having the following properties:
\begin{enumerate}
\item
Selfadjoint: $\rho=\rho^{*}$, where $\rho^{*}$ is the adjoint of $\rho$.
\item
Positive: $\rho\geq 0$, or equivalently $\langle\psi, \rho\psi\rangle\geq 0$ for all $\psi\in\mathcal{H}$.
\item
Trace one: $\Tr(\rho)=1$.
\end{enumerate}
Positivity implies that the eigenvalues of $\rho$ are all nonegative and by the last
property, they sum up to one. Notice that the above requirements parallel those of
defining probability densities.

We have the following diagonal form
\begin{equation}
\rho=\sum_{i=1}^{{\rm dim}\mathcal{H}}\lambda_{i}\rho_{i}
\end{equation}
where $\rho_{i}$ is the projection onto the one dimensional space generated by the
eigenvector $e_{i}\in\mathcal{H}$ of $\rho$ corresponding to the eigenvalue $\lambda_{i}$, i.e.,
$\rho e_{i}=\lambda_{i}e_{i}$. With respect to a fixed orthonormal basis
$\{\psi_{i}\}_{i\geq 1}$ in $\mathcal{H}$, the operator $\rho$ can be represented as a
matrix with elements $\rho_{i,j}=\langle\psi_{i}, \rho\psi_{j}\rangle$.


Let us consider now the following problem.
We are given a quantum system prepared in an unknown state $\rho$ and we would like
to determine $\rho$. In order to obtain information about the system we have to
measure its properties. The laws of quantum mechanics say that for any 
given measurement with space of outcomes given by the measure space
$(\mathcal{X}, \Sigma)$, the result of the measurement performed on a system
prepared in state $\rho$ is {\it random} and has probability distribution
$\mathbb P_{\rho}$ over $(\Omega, \Sigma)$ such that the map
\begin{equation}
\rho\mapsto\mathbb P_{\rho},
\end{equation}
is affine, i.e. it maps convex combinations of states into the corresponding convex
combination of probability distributions. This has a natural interpretation:
a system can be prepared in a mixture $\lambda\rho_{1}+(1-\lambda)\rho_{2}$
of states by randomly choosing the preparation procedure
according to the individual state $\rho_{1}$ with probability $\lambda$ and $\rho_{2}$
with probability $1-\lambda$. The distribution of the results will then reflect this
randomized preparation as well.

The most common measurement is that of an observable such as
energy, position, spin, etc.
Any given observables is described by some {\it selfadjoint operator} $\mathbf{X}$
on the Hilbert space $\mathcal{H}$ and we suppose for simplicity that it has a discrete
spectrum, that is, it can be written in the diagonal form
\begin{equation}
\mathbf{X}=\sum_{i=1}^{{\rm dim}\mathcal{H}}x_{a} \mathbf{P}_{a}.
\end{equation}
with $x_{a}\in\R$ and $\mathbf{P}_{a}$ one dimensional projections onto the eigenvectors of
$\mathbf{X}$. The result $X$ of the measurement of the observable $\mathbf{X}$ for a
preparation given by the state $\rho$, is a random element of the set
$\Omega=\{x_{1}, x_{2}, \dots\} $ of eigenvalues of $\mathbf{X}$ and has the probability distribution
\begin{equation}
\mathbb{P}_{\rho}\left[  X=x_{a} \right]=\Tr( \mathbf{P}_{a}\rho).
\end{equation}
This measurement will give statistical information about the diagonal elements of the density matrix $\rho$ with respect to the eigenbasis of ${\bf X}$, and it suggests 
that in order to estimate all matrix elements of $\rho$ one would have to probe the system from a number of ``directions''  by performing different measurements on identically prepared systems. 
 
This can be generalized to the case of infinite dimensional Hilbert space, and measurements with outcomes in arbitrary measure spaces. In the next section we will see that an infinite dimensional density matrix can be estimated by measuring a randomly chosen observable from a continuous family of non-commuting observables.


 \subsection{Quantum homodyne tomography and the Wigner function}
\label{sec.qhomodyne}

An important example of quantum system is the monochromatic light in a cavity, described
by density matrices on the Hilbert space of complex valued square integrable
functions on the real line, $L_{2}(\R)$. A distinguished orthonormal basis of this space is given by the vectors
\begin{equation}\label{eq.psi_n}
\psi_k(x)=H_k(x)e^{-x^2/2},\qquad k =0,1, 2\dots,
\end{equation}
where $H_{k}$ are the Hermite polynomials normalized such that $\psi_k(x)$ is a unit
vector representing the pure state of $k$ photons inside the cavity.
We will denote the matrix elements of $\rho$ with
respect to this basis by $\rho_{i,j}$. Notice that the diagonal of the density matrix
is a probability distribution over the nonnegative numbers $p_{k}=\rho_{k,k}$.
This is the probability distribution of results when measuring the number of photons in the cavity prepared in state $\rho$. Clearly this distribution does not contain information about the off-diagonal elements of $\rho$ thus it is not sufficient for identifying the state
of the system. This is a typical situation in state estimation and one has to devise
experiments in which the systems are looked at subsequently from ``different directions'',
a broadly described methodology which in the physics literature goes by the name of
quantum tomography. Quantum Homodyne Tomography is one such measurement method which is frequently used in quantum optics at the estimation of the quantum state of light \cite{Breitenbach&Schiller&Mlynek,Ourjoumtsev,Zavatta}. 
We skip the measurement set-up which is described in detail in \cite{Artiles&Gill&Guta,Butucea&Guta&Artiles} and present the statistical problem 
associated to this measurement.

We observe $(X_1,\Phi_1),\dots ,(X_n,\Phi_n)$, i.i.d.~random variables with values in
$\mathbb{R}\times [0,\pi]$ and  distribution $\mathbb P_\rho$ whose density with respect to the  measure  $\pi^{-1}d\phi \times dx$ is given by 
\begin{equation}\label{eq.concrete.p}
p_\rho(x,\phi)=\sum_{j,k=0}^\infty \rho_{j,k}\psi_k(x)\psi_j(x)e^{-i(j-k)\phi}.
\end{equation}
Since $\rho$ is a positive definite matrix of trace $1$, and $\psi_j$ form an orthonormal basis, it follows that $p_{\rho}$ is a probability density: real, nonnegative, integrates to $1$. The data $(X_\ell,\Phi_\ell)$, $\ell=1,\ldots,n$, come from independent measurements on identically prepared pulses of light escaping from the cavity, whose state is completely
encoded in the matrix $\rho$.

For each of the systems independently, we repeat the following experimental 
procedure: we first choose the angle $\Phi$ uniformly distributed over $[0,\pi]$ and then measure a certain observable ${\bf X}_{\phi}$ called quadrature, obtaining a real valued result with probability density $p_{\rho}(x,\phi)$. The quadrature is defined as the linear combination $X_{\phi}:= \cos \phi {\bf Q} + \sin \phi {\bf P}$, where ${\bf Q}$ and ${\bf P}$ are the electric and magnetic fields of the light beam given by the selfadjoint operators
\begin{eqnarray*}
&& \mathbf{Q}\psi(x)=x\psi(x), \qquad {\rm and } \qquad\mathbf{P}\psi(x)=-i\frac{d\psi}{dx}.
\end{eqnarray*}
The characteristic functions of these densities can be put together to define a function of two variables
\begin{equation}\label{eq.Fourier transform.tildeW}
\widetilde{W}_\rho(u,v) :=
\mathrm{Tr}\big(\rho \exp(-it\mathbf{X}_\phi)\big)=
\mathcal{F}_1[p_\rho(\cdot,\phi)](t),
\end{equation}
where we have used the polar coordinates $(u,v)=(t\cos\phi, t\sin\phi)$, and  
$\mathcal{F}_1$ is the Fourier transform  with respect to the first variable, for fixed 
$\phi$. Note that our convention for defining the Fourier transform and its inverse are the following 
$$
\mathcal{F}[f] (t) = \int_{-\infty}^{\infty} f(x) e^{-ixt}dx,
\qquad
\mathcal{F}^{-1}[g] (x) = \frac{1}{2\pi} \int_{-\infty}^{\infty} g(t) e^{ixt}dt.
$$
Equivalently, we can write 
\begin{equation}\label{def.Wigner}
\widetilde{W}_\rho(u,v)=
\mathrm{Tr}\big(\rho \exp(-iu\mathbf{Q}-iv\mathbf{P})\big),
\end{equation}
which resembles a characteristic function of a bivariate probability density, namely the joint density of ${\bf Q}$ and ${\bf P}$. However, since the operators 
${\bf Q}$ and ${\bf P}$ do not commute with each other, we cannot speak of their joint probability distribution and the function $\widetilde{W}_\rho(u,v)$ is in general not the characteristic function of a probability density but rather of the so called 
Wigner function 
\begin{equation}\label{eq.Wigner.fct}
W_{\rho} (q,p):=\mathcal{F}_2^{-1}[\widetilde W_{\rho}](q,p),
\end{equation}
a ``quasi-distribution'' which may take negative values but whose marginals are bona-fide probability densities. As we will see below, the Wigner function $W_{\rho}$ is in one to one correspondence with the density matrix $\rho$, and in quantum optics one frequently uses the Wigner function as an alternative representation of the 
quantum state, having the advantage that it illustrates important  ``quantum features''  such as squeezing and negative oscillations. 
From \eqref{eq.Fourier transform.tildeW} and \eqref{eq.Wigner.fct} we deduce that 
the probability density of the data $p_{\rho}(x,\phi)$ is the {\it Radon transform} of the Wigner function 
\begin{equation*}
\mathcal{R}[W_{\rho}](x,\phi)=\int_{-\infty}^\infty
W_{\rho}( x\cos\phi - t\sin\phi,  x \sin\phi + t\cos\phi )dt.
\end{equation*}
adding Quantum Homodyne Tomography to a long the list of applications 
ranging from computerized tomography to astronomy and geophysics \cite{Deans}.

Another important feature of the Wigner function is that it can be used as a computational tool:  for any selfadjoint operator $\mathbf{X}$ there exists a function 
$W_\mathbf{X}$ from $\R^2$ to $\R$ such that the expectation of ${\bf X}$ is given by
\begin{equation}
\mathrm{Tr}(\mathbf{X}\rho)=2\pi \iint W_\mathbf{X}(q,p)W_\rho(q,p)dqdp.
\end{equation}
In particular, the correspondence between the density matrix $\rho$ and the Wigner function $W_{\rho} $ is an $L_{2}$ isometry up to a constant
\begin{equation}\label{eq.isometry}
\| W_{\rho}- W_{\tau}\|_2^2 :=\iint |W_{\rho}(q,p)- W_{\tau}(q,p)|^2dp dq=
\frac{1}{2\pi}\sum_{i,j=0}^{\infty} |\rho_{i,j}-\tau_{i,j}|^2.
\end{equation}


The space of Wigner function has an overlap with that of probability distributions. For example, all Gaussian densities which are bounded from above by $1/\pi$ are Wigner functions, while the rest do not correspond to physical states. This is due to the celebrated {\it Heisenberg's uncertainty relations} which say that the non-commuting observables $\mathbf{P}$ and $\mathbf{Q}$ cannot have probability
distributions such that the product of their variances is smaller than $\frac{1}{4}$. In general, a Wigner function cannot be too ``peaked'':
\begin{equation}\label{eq.uniform.bounded.Wigner}
| W_{\rho}(q,p)|  \leq \frac{1}{\pi}, \qquad \mathrm{for~all}~(q,p)\in\R^2.
\end{equation}

Some examples of quantum states which can be created in laboratory are given in
Table 1 of \cite{Artiles&Gill&Guta}. Typically,  the corresponding Wigner functions
have a Gaussian tail but need not be positive. For example the state of one-photon in the cavity
is described by the density matrix with $\rho_{1,1}=1$ and all other elements zero which is equal
to the orthogonal projection onto the vector $\psi_{1}$. The corresponding Wigner function is
\begin{equation*}
W(q,p)=\frac{1}{\pi} (2q^2+2p^2-1) \exp (-q^2-p^2).
\end{equation*}

In conclusion, although we deal with a problem which is similar to that of Positron Emission Tomography, the parameter space is different from the space of probability densities and special techniques have to be developed for this situation. 

\section{The main results}
\label{sec.stat.setup}

Our problem is that of estimating the Wigner function $W_\rho(z)$, defined on the plane $z=(q,p)$. In order to prove rates of convergence some restrictions are necessary to be imposed to the class.
We consider the class $\cW(\beta, L)$ of Wigner functions which are continuous and whose Fourier transform satisfy
\beq
\frac{1}{4\pi^{2}}\int_{\R^{2}} |\widetilde W (w)|^{2} \exp(2\beta |w|) \, dw \leq L
\eeq
for $\beta$ and $L$ positive constants. This condition implies  that the function to be estimated is very smooth.
Such classes appeared in the statistical literature in \cite{Ibragimov&Hashminski.1}, and we subsequently used in the context of density estimation \cite{Golubev.Levit.}, 
functional estimation, regression problems \cite{Golubev.Levit.Tsybakov}, and  tomography \cite{Cavalier}. In \cite{Butucea&Guta&Artiles} we have argued that from the physical point of view it is natural to assume that the Wigner function of a state which can be created in the lab belongs to such a class.

We use a kernel-type estimator based on the following function called a 
band-limited filter 
\beq
K_{\delta_{n}}(u)=\frac{1}{4\pi}\int_{-\delta_{n}}^{\delta_{n}} r e^{iru}.
\eeq
This filter has already been used in the context of tomography \cite{Natterer,Korostelev&Tsybakov,Cavalier}, and its Fourier transform is
\beq \label{prop.kernel}
\widetilde K_{\delta_{n}} (t)=\frac{1}{2}|t| I_{\delta_{n}}(t),
\eeq
where $I_{\delta_{n}}$ is the indicator function of $\{t: |t|\leq 1/\delta_{n}\}$.
The estimator we use is
\beq
\widehat W_{n}(z)= \frac{1}{n}\sum_{i=1}^{n} K_{\delta_{n}}([z,\Phi_{i}]-X_{i}  )
\eeq
with i.i.d. observations $(X_{i}, \Phi_{i})$, for $i=1,\dots, n$, with density
$p_{\rho}(x,\phi)=
\cR [W_{\rho}](x,\phi)$. \\
Following \cite{Natterer} we define the dual operator $\cR^{\#}$ on $L_{1}(\R\times [0,\pi])$ by
\beq
\mathcal{R}^{\#}[h](z)=\int_{0}^{2\pi} h([z,\phi], \phi) \, d\phi.
\eeq
Then
\beq
\cR^{\#}\cR [W](z)=\int_{0}^{2\pi} \cR[W] ([z,\phi], \phi) \,d\phi
\eeq
represents the integrals of $W$ over all lines passing through the point $x$.
Note that in general $R^{\#}\cR [W](z)\geq 0$ for all Wigner functions $W$ and all $z\in\R^{2}$, and the number states $\psi_{k}$ with $k$ odd have the property that
$\cR^{\#}\cR [W](0)=0$. In \cite{Cavalier} it is assumed that the probability distributions $f$ to be estimated are strictly positive which implies that $\cR^{\#}\cR [f](z)> 0$ for all 
$z$. For the upper bound we will assume that the latter condition holds.


\subsection{The upper Bound}

\begin{theorem}\label{theorem.upperbound}
For any $W\in\cW (\beta, L)$ and any fixed $z\in\R^{2}$ such that $R^{\#}\cR [W](z)>0$ we have as $n\to\infty$,
\beq
\bE \left[ (\widehat{W}_{n}(z)-W(z))^{2} \right]=
C^{*} \cR^{\#}\cR[W](z)\times \frac{(\log n)^{3}}{n}(1+o(1))
\eeq
where $C^{*}= \frac{\pi}{3 (4\pi\beta)^{3}}$.
\end{theorem}

\noindent
\textit{Proof.} 
We will provide only the main steps of the proofs pointing out where the assumption on the class of Wigner function plays a role. For a more detailed proof of the bounds for
the class of probability densities in $\cA(\beta, L)$ we refer to \cite{Cavalier}. The risk can be decomposed in two parts, the bias and the stochastic part
\begin{align}
\bE \left[ (\widehat{W}_{n}(z)-W(z) )^{2} \right] &=
\left( \bE (\widehat{W}_{n}(z)) - W(z) \right)^{2} +
\bE \left[  (\widehat{W}_{n}(z)- \bE(\widehat{W}_{n}(z)) )^{2}\right] \nonumber \\
& :=
b_{n}^{2}(z)+\sigma_{n}^{2}(z).
\end{align}
On one hand, using property (\ref{prop.kernel}) of the kernel and the inverse Fourier transform,  the bias can be written as
\beq
b_{n}(z)=\frac{1}{(2\pi)^{2}} \int \widetilde{W}(w) I\left(|w|>\frac{1}{\delta_{n}}\right)
e^{-i\br w, z \ke} \, dw.
\eeq
By using Cauchy-Schwarz and the fact that $W\in\cW(\beta,L)$ we get
\beq
b_{n}^{2}(z)\leq \frac{1}{\delta_{n}} e^{-2\beta /\delta_{n}}(1+o(1)),
\eeq
as $\delta_{n}\to 0$. With the choice $1/\delta_{n}=\log n/(2\beta)$ the bias upper bound becomes
\beq
b_{n}^{2}(z)\leq c\frac{\log n}{n}(1+o(1)), \qquad \mathrm{as}~n\to \infty.
\eeq
On the other hand, the variance is equal to
\begin{align}
\sigma_{n}^{2}(z) = & \frac{1}{n}{\rm Var} K_{\delta_{n}}([z,\Phi]-X )=
\nonumber\\
&
\frac{1}{n}\bE \left[ K_{\delta_{n}}^{2}([z,\Phi]-X)\right]-
\frac{1}{n}\left( \bE \left[ K_{\delta_{n}}([z,\Phi]-X)\right]\right)^{2},
\end{align}
where $(X,\Phi)$ is a random variable with probability density $p_{\rho}(x,\phi)=
\cR[W_{\rho}](x,\phi)$.
The second term can be bounded as follows
\beq
\frac{1}{n}\left( \bE \left[ K_{\delta_{n}}([z,\Phi]-X)\right]\right)^{2}\leq
\frac{1}{(2\pi)^{2}n} \int  I_{\delta_{n}}(|w|)e^{-2\beta|w|} \, dw
\times \int |\widetilde W (w)|^{2} e^{2\beta |w|} \, dw = O(\frac{1}{n}),
\eeq
with $O(\frac{1}{n})$ uniformly with respect to $W\in\cW(\beta,L)$.\\
The first term is
\begin{equation}
\bE\left[K_{\delta_{n}}^{2}([z,\Phi]-X)\right]=
\int_{0}^{\pi} \int_{\R} K_{\delta_{n}}^{2}([z,\phi]-y)p_{\rho}(y,\phi)  \, d\phi \, dy.
\end{equation}
Now denote
\beq
G(u)=\frac{3}{\pi}\left( \int_{0}^{1} r \cos(ur) \, dr\right)^{2},
\eeq
and let $G_{\delta}(u)=(1/\delta) G(u/\delta)$. We have
\beq
\int_{\R} K_{\delta_{n}}^{2}([z,\phi]-y)p_{\rho}(y,\phi)  \, dy=
\frac{\pi^{2}}{3(2\pi)^{4}} \left(\frac{1}{\delta_{n}}\right)^{3}
( G_{\delta_{n}}\ast \cR[W](\cdot, \phi) )([z,\phi]).
\eeq
Using \cite{Cavalier} Lemma 4, we have that as $\delta\to 0$,
\bec
\int_{0}^{\pi} (G_{\delta} \ast  \cR[W](\cdot, \phi) ) \, d\phi =
\cR^{\#}\cR [W](z) (1+o(1))+O(\delta^{1/3}).
\eec
With this the second term of the variance can be written as
\bec
\bE\left[K_{\delta_{n}}^{2}([z,\Phi]-X)\right]=
\frac{\pi^{2}}{3(2\pi)^{4}} \left(\frac{1}{\delta_{n}}\right)^{3}\cR^{\#}\cR [W](z) (1+o(1))
+O\left(\left(\frac{1}{\delta_{n}}\right)^{4-4/3}\right),
\eec
as $\delta_{n}\to 0$. Thus as $n\to \infty$
\begin{align}\label{eq.variance}
\sigma_{n}^{2}(z) & \leq C^{*}\cR^{\#}\cR [W](z) \frac{(\log n)^{3}}{n}(1+o(1))+
O\left( \frac{(\log n)^{4-4/3}}{n}\right) \nonumber \\
& = \frac{\log {}^3 n }{n} \left( C^{*}\cR^{\#}\cR [W](z) (1+o(1))+
O(\log ^{-1/3} n) \right) \nonumber
\end{align}
\qed

\noindent If $R^{\#}\cR [W](z)>0$ then we obtain the claimed constant. Notice that if $R^{\#}\cR [W](z)=0$ then the convergence rate is faster than $\frac{\log {}^3 n }{n}$. \\

\subsection{The lower bound}

In order to prove a lower bound result we consider the slightly modified class of Wigner functions
\beq
\cW(\beta, L,\alpha_{n})=\{ W \in\cW(\beta,L): ~ \cR^{\#}\cR [W](z)\geq \alpha_{n}\} ,
\eeq
for a sequence $\alpha_{n}$ such that $\lim_{n\to\infty}\alpha_{n}=0$ and
$\lim_{n\to \infty}(\alpha_{n}(\log n)^{1/3})=\infty$. Let us denote
\beq \label{def.rate}
r_{n}(W, z)=
\left(
C^{*}\cR^{\#}\cR [W](z)\frac{(\log n)^{3}}{n}
\right)^{1/2}.
\eeq
\begin{theorem}\label{theorem.lowerbound}
For a fixed $z\in\R^{2}$, we have
\bec
\liminf_{n\to \infty} ~\inf_{\widehat{W}_{n}} ~\sup_{W\in\cW(\beta,L,\alpha_{n})}
\bE \left[ \left(
\frac{\widehat{W}_{n}(z)-W(z)}{r_{n} (W, z)}
\right)^{2}\right]
\geq 1
\eec
where $\inf_{\widehat{W}_{n}}$ denotes the infimum over all estimators of $W(z)$.
\end{theorem}

\noindent
\textit{Proof.} 
The proof is based on the standard procedure of building a hardest parametric 
subfamily for the class  $\cW(\beta, L)$  of Wigner functions of the form
\beq
W_{c}=W_{\alpha}+cg_{a},
\eeq
where $c$ is a parameter in a neighborhood of the origin, $W_{\alpha}$ and 
$g_{a}$ are functions to be defined shortly.

The essential point of the proof is that the family of probability densities 
$f_{0}+ cg_{a}$ used in \cite{Cavalier} is not always contained in our parameter space consisting of Wigner functions. For illustration we will show that for some parameters $\beta$, the function $f_{0}$ defined in equation (29) of \cite{Cavalier} is not the Wigner function of a quantum state. 
Indeed, suppose for the moment that this was the case, i.e. 
$f_{0} = W_{\rho}$ for some state $\rho$. Then by the rotation symmetry of $f_{0}$, 
the density matrix $\rho$ must be diagonal and  
\bec
p_{\rho}(x,\phi)=  \sum_{k=0}^{\infty} \rho_{k,k}\psi_{k}^{2}(x).
\eec
Using the inequality \cite{Erdelyi} $\|\psi_{k}\|_{\infty}\leq k$, where $k$ is a constant whose value is slightly bigger that $1$, and the fact that $\sum \rho_{k,k}=1$, we find that $\| p_{\rho}\|_{\infty} \leq k$. However, the Radon transform of the Wigner function 
$W_{\rho}$ is 
\bec
p_{\rho} (x,\phi) = \cR[W_{\rho}](x,\phi) =\cR[f_{0}](x,\phi)=\frac{\beta}{\pi} \frac{1}{x^{2}+\beta^{2}}, 
\eec
which violates the above bound for $\beta\leq 1/(\pi k)$.\\

We thus define a parametric subfamily of $\cW(\beta, L)$ which is a suitable modification of the family considered in \cite{Cavalier}  in order to cope with this
problem.\\

\noindent {\it Construction of $W_{\alpha}$.}
Consider the Mehler formula, (see \cite{Erdelyi}, 10.13.22)
\begin{equation}\label{eq.Mehler}
\sum_{n=0}^\infty z^n \frac{1}{\sqrt{\pi}n!2^n} H_n(x)^2 e^{-x^2}=
\frac{1}{\sqrt{\pi(1-z^2)}}\exp\left(-x^2\frac{1-z}{1+z}\right).
\end{equation}
Integrating both terms with $f_{\alpha}(z)=\alpha(1-z)^{\alpha}$ we get
\begin{equation}
p_{\alpha}(x,\phi):=\sum_{n=0}^\infty   \psi_n (x)^2 \int_{0}^{1}f_{\alpha}(z)z^{n} \, dz =
\int_{0}^{1}\frac{f_{\alpha}(z)}{\sqrt{\pi(1-z^{2})}} \exp\left(-x^{2}\frac{1-z}{1+z}\right) \, dz.
\end{equation}
The Fourier transform of $p_{\alpha}$ is
\beq\label{eq.W0}
\widetilde{W}_{\alpha}(w):=\mathcal{F} [p_{\alpha}](w)=
\int_{0}^{1} \frac{f_{\alpha}(z)}{1-z}\exp\left(-|w|^{2}\frac{1+z}{4(1-z)}\right) \, dz
\eeq
Notice that the normalization condition $\int p_{\alpha}=1$ is equivalent to $\widetilde{W}_{\alpha}(0)=1$ which is satisfied for the chosen functions $f_{\alpha}$, thus $p_{\alpha}$ is a probability density corresponding to a diagonal density
matrix $\rho^{\alpha}$ with elements
\beq
\rho^{\alpha}_{k,k}=\int_{0}^{1}z^{k}f_{\alpha} \, dz.
\eeq
We denote by $W_{\alpha}$ the Wigner function whose Fourier transform is defined in equation (\ref{eq.W0}) with $\alpha>0$ a parameter to be fixed later. This function is considered in a more general form in \cite{Butucea&Guta&Artiles}, and corresponds 
to $f_{\alpha}^{\epsilon}$ for  $\epsilon = 0$ .

\noindent {\it Construction of $g_{a}$.} Let \cite{Cavalier}
\beq \label{eq.Ha}
H_{a}(s)=\frac{1}{(2\pi)^{2}}\int_{0}^{\infty}\frac{r}{1+a^{-1}\sinh^{2} \beta r}\cos(sr) \, dr,
\qquad s\in\R,
\eeq
where $a>0$ is a parameter which will depend on $n$ as $a=a_{n}=n^{\eta}$ with $0<\eta <1$. 
The Fourier transform of $g_{a}$ is
\beq
\widetilde{H}_{a}(t)=\frac{1}{4\pi}\frac{t}{1+a^{-1}\sinh^{2} \beta t}, \qquad t \in\R.
\eeq
Let
\beq \label{eq.g_a}
g_{a}(z)=\frac{1}{2(2\pi)^{3}}\int_{\R^{2}}\frac{|w|}{1+a^{-1}\sinh^{2} \beta |w|}
\cos(\br z, w\ke) \, dw, \qquad z\in\R^{2},
\eeq
and its Fourier transform
\beq
\widetilde{g}_{a}(w)=\frac{1}{4\pi}\frac{|w|}{1+a^{-1}\sinh^{2} \beta |w|},
 \qquad w \in\R^{2}.
\eeq
%

Now, let us consider the family
\beq
W_{c}=W_{\alpha}+cg_{a}
\eeq
where the real parameter $c$ satisfies
\beq
|c|\leq C_{a}=\frac{q}{\sqrt{a}(\log a)^{3/2}},
\eeq
with  $q>0$ sufficiently small.
By translating with $z$ in $\R^{2}$ we obtain our hardest family for pointwise estimation at the point $z$
\beq \label{eq.shift}
W_{c}^{z}(\zeta)=W_{c}(\zeta-z).
\eeq
We will check now that $W_{c}^{z}$ belongs to the class $\cW(\beta, L)$ for an appropriate choice of $\alpha$ in  $W_{\alpha}$, which means that $W_{c}^{z}$ is 
a Wigner function and
\beq\label{eq.norm}
\frac{1}{4\pi^{2}}\int \left|\widetilde W_{c}^{z} (w) \right|^{2} e^{2\beta|w|} \, dw = 
\frac{1}{4\pi^{2}}\int \left|\widetilde W_{c} (w) \right|^{2} e^{2\beta|w|} \, dw \leq L.
\eeq

By Lemma \ref{lemma.W0class} we have $W_{\alpha}\in \cW(\beta,L/4)$  for a small enough $\alpha>0$, and from Lemma 5 of \cite{Cavalier} we have that
\beq
\frac{c^{2}}{4\pi^{2}}\int \left| \tilde g_{a} (w) \right|^{2} e^{2\beta|w|} \, dw \leq L/4,
\eeq
for all $c\leq C_{a}$. The last two conditions together imply (\ref{eq.norm}).

Furthermore, $W_{c}^{z}$ has to be a Wigner function. As translations in the plane transform Wigner functions into Wigner functions, we need only to show this for $W_{c}$. This means that there exist a family of density matrices $\rho_{c}$ such that their corresponding Wigner functions are $W_{c}$. The invariance of $W_{c}$ under rotations in the plane translates into the fact that $\rho^{c}$ has all off-diagonal elements equal to zero, and thus we only need to show that all its diagonal elements are positive and add up to one.  The relation between the diagonal matrix elements and the Wigner function 
is \cite{Leonhardt}
\beq
\rho^{c}_{k,k} = \frac{1}{2\pi} 
\int_{\R^{2}} e^{|w|^{2}/4} L_{k}(|w|^{2}/2)\widetilde{W}_{c}(w) \, dw,
\eeq
where $L_{k}$ are the Laguerre polynomials. By linearity we have
\beq
\rho^{c}_{k,k}=\rho^{\alpha}_{k,k}+c\tau^{a}_{k,k}
\eeq
where
\beq\label{eq.tau}
\tau^{a}_{k,k}= 
\int_{0}^{\infty}t e^{-t^{2}/4} L_{k}(t^{2}/2)\widetilde{H}_{a}(t) \, dt,
\eeq
and
\beq
\rho^{\alpha}_{k,k} = \alpha\int_{0}^{1} z^{k}(1-z)^{\alpha}.
\eeq
Corroborating the result shown in Lemma \ref{lemma.rhoa}
\bec
\tau^{a}_{kk} =O\left( k^{-5/4}(\log a)^{4} \right),
\eec
as $a$, $k\to\infty$, with that of Lemma 2 in \cite{Butucea&Guta&Artiles}
\bec
\rho^{\alpha}_{k,k}\sim k^{-(1+\alpha)},
\eec
we conclude that if $\alpha<1/4$, then $\rho^{c}_{k,k}\geq 0$ for all $k\geq 0$ and
$|c|\leq C_{a}$ for $a$ sufficiently large.

\noindent Now we can use the fact that for the family of translated functions, as defined in eq. (\ref{eq.shift}), $\cR^{\#}\cR[W^z](z)=\cR^{\#}\cR[W](0)$. Indeed, using \eqref{eq.Fourier transform.tildeW} we get
\begin{eqnarray}
\cR [W^z](x,\phi)&&= \frac{1}{2\pi} \int \widetilde{W}^z(t\cos\phi, t\sin\phi) e^{itx} \,dt \nonumber\\
&&= \frac{1}{2\pi} \int \widetilde{W}(t\cos\phi, t\sin\phi) e^{-it[z,\phi]} e^{itx} \,dt \nonumber\\
&&= \frac{1}{2\pi} \int \widetilde{W}(t\cos\phi, t\sin\phi) e^{it(x-[z,\phi])} \,dt \nonumber\\[2ex]
&&= \cR [W](x-[z,\phi],\phi). \label{eq.translated}
\end{eqnarray}
Now, from definition of $\cR^{\#}$ we get
\begin{eqnarray}
\cR^{\#}\cR[W^z](z)
&=& \int_{0}^{2\pi} \cR[W^z]([z,\phi], \phi)d\phi \nonumber\\
& = & \int_{0}^{2\pi} \cR [W](0,\phi) d\phi \ = \ \cR^{\#}\cR[W](0) \nonumber.
\end{eqnarray}
Thus, $R_c(z) : = \cR^{\#}\cR[W_c^z](z) = \cR^{\#}\cR[W_\alpha](0) + c \cR^{\#}\cR[g_a](0)$ for any $z$. Given in our case $W_\alpha$ and $g_a$ are invariant under rotations we obtain
\begin{eqnarray}
R_0 & := & \cR^{\#}\cR[W_\alpha](0) = \int_{0}^{2\pi} \cR [W_\alpha](0,\phi) \, d\phi = \int_{0}^{2\pi} p_\alpha(0,\phi) \, d\phi \nonumber \\
& = & 2\sqrt{\pi} \alpha\int_0^1 \frac{(1-z)^{\alpha-1/2}}{(1+z)^{1/2}} \, dz > 0.
\label{eq.r0}
\end{eqnarray}
For the second term, using eq. (\ref{eq.Ha}), Lemma 5 in \cite{Cavalier}, and definition of $C_a$
\beq \label{eq.sup}
\sup_{|c|\leq C_a} | c \cR^{\#}\cR[g_a](0)| = \sup_{|c|\leq C_a} |2 \pi c H_a(0) | = o(1) \quad \text{as} \quad a \to \infty.
\eeq
We conclude that $R_c = R_0(1+o(1)) \geq \alpha_n,$ for $n \to \infty$ and thus, for $a$ large enough, $W_c^z \in \cW(\beta,L,\alpha_n)$. Moreover, from (\ref{def.rate}), 
\begin{equation}\label{eq.rn.r0}
r_{n}(W_c^z,z)^{2} = r_{n}(W_\alpha,0)^{2}(1+o(1)) = C^{*}R_{0}\frac{(\log n)^{3}}{n}(1+o(1)).
\end{equation}
The rest of the proof is based on the Van Trees inequality and follows along the lines of \cite{Cavalier}. The main difference is in the proof of Lemma \ref{lemma.information} where the Fisher information of the family of densities $\frac{1}{\pi} \cR [W_c^z]$ is approximated. Take a continuously differentiable probability density, $\lambda_0(c)$, defined on the interval $[-1,1]$, such that $\lambda_0(-1)=\lambda_0(1)=0$, with a finite Fisher information $I_0$. The new density $\lambda(c)=\lambda_a(c) = C_{a}^{-1}\lambda_0(C_{a}^{-1}c)$ is a prior density with finite Fisher information $I(\lambda)=I_0 C_a^{-2}$. Finally let us  define $I(c)$ the Fisher information of the family of densities $\frac{1}{\pi}\cR [W_c^z]$. Using the Van Trees inequality
\begin{eqnarray}
\inf_{\widehat{W}_{n}}  ~\sup_{W\in\cW(\beta,L,\alpha_{n})}
\bE_{W} \left[ \left(
\widehat{W}_{n}(z)-W(z)
\right)^{2}\right] & \geq &
\inf_{\widehat{W}_{n}} ~\sup_{|c|<C_{a}}
\bE_{W_c^z} \left[ \left(
\widehat{W}_{n}(z)-W_c^z(z)
\right)^{2}\right] \nonumber \\
 \geq \inf_{\widehat{W}_{n}} \int_{-C_a}^{C_a}
\bE_{W_c^z} \left[ \left(
\widehat{W}_{n}(z)-W_c^z(z)
\right)^{2}\right] \lambda_{a}(c)dc &\geq& 
\frac{(\int_{-C_a}^{C_a} 
\left(
\partial W_c^z(z)/ \partial c)\lambda_{a}(c)\, dc
\right)^2}
{n \int_{-C_a}^{C_a} I(c)\lambda_{a}(c)\, dc + I(\lambda)}.
\end{eqnarray}
Now, from equation (\ref{eq.g_a}) and Lemma 5 of \cite{Cavalier}, we get that $\partial W_c^z(z)/ \partial c = g_a(0) = C^{*} (\log a)^3(1+o(1))$. By inserting the expression of 
$I(c)$ from Lemma \ref{lemma.information}, and using \eqref{eq.rn.r0} and 
$a=a_{n}= n^{\eta}$ we get the lower bound 
\begin{eqnarray*}
 \inf_{\widehat{W}_{n}} ~\sup_{|c|<C_{a}} \bE \left[ \left( \frac{\widehat{W}_{n}(z)-W^{z}_{c}(z)}{r(W_c^z,z)} \right)^{2}\right] 
 &&\geq 
 \frac{(\int_{-C_a}^{C_a} (\partial W_c^z(z)/ \partial c)\lambda_c\, dc)^2}{n \int_{-C_a}^{C_a} I(c)\lambda_c\, dc + I(\lambda)} \cdot 
 (r_{n}(W_c^z,z)(1+o(1)))^{-2} \\
 &&= 
 \frac{(C^{*} \eta^3 (\log n)^3  )^{2}\, r_n(W_\alpha,0)^{-2}(1+o(1) )}{C^{*} \eta^3 R_0^{-1} n (\log n)^3 (1+o(1))+ q^{-2} I_0 \eta^{3} n^\eta (\log n)^3}. \nonumber
\end{eqnarray*}
By letting $n\to \infty$ followed by  $\eta\to 1$ we finally obtain
\begin{eqnarray*}
\liminf_{n\to\infty}\, \inf_{\widehat{W}_{n}} ~\sup_{W\in\cW(\beta,L,\alpha_{n})} \bE \left[ \left( \frac{\widehat{W}_{n}(z)-W(z)}{r(W,z)} \right)^{2}\right] &&\geq 1. \nonumber
\end{eqnarray*}

\qed

\section{Technical Lemmas}\label{sec.lemmas}

\begin{lemma}\label{lemma.w}
Let $W\in \mathcal{W}(\beta,L)$. Then the following inequalities hold
\begin{eqnarray}
&&
|\widetilde{W}(\omega)|\leq 1, \qquad \forall \omega\in \R^{2},\\
&&
\int_{\R^{2}}|\widetilde{W}(w)| \, dw \leq Q,\label{eq.second}\\
&&
\sup_{W \in\mathcal{W}(\beta,L)} \cR^{\#}\cR[W](z)\leq 1+\frac{Q}{2\pi},
\end{eqnarray}
where $Q$ is a constant depending only on $\beta$ and $L$.
\end{lemma}

\noindent\textit{Proof.} As $W=W_{\rho}$ for some density matrix $\rho$ we use the definition of
$\widetilde{W}_{\rho}$ to obtain
\bec
|\widetilde{W}_{\rho}(u,v)|=|\Tr \left( \rho e^{-iu\mathbf{Q}-iv\mathbf{P}}\right)| \leq ||\rho||_{1}~||e^{-iu\mathbf{Q}-iv\mathbf{P}}||=1,
\eec
where we have used the normalization of the density matrix $\rho$ and the fact that $\|U\|=1$ for any unitary operator $U$. The inequality \eqref{eq.second} is a direct consequence of the definition of
$\mathcal{W}(\beta, L)$ and the Cauchy-Schwarz inequality.
Now using (\ref{eq.Fourier transform.tildeW}) and the previous inequalities we get
\begin{eqnarray}
|\cR^{\#}\cR[W](z)|&&=
\left|\int_{0}^{2\pi} \cR[W]([z,\phi], \phi)d\phi\right|=
\nonumber\\
&&
\left|\frac{1}{2\pi}\int_{0}^{2\pi}d\phi \int dt \widetilde{W}(t\cos\phi, t\sin\phi) e^{it[z,\phi]} \right|\leq
\nonumber\\
&&
\frac{1}{2\pi} \int_{0}^{2\pi}d\phi \int dt |\widetilde{W}(t\cos\phi, t\sin\phi)|\leq
\nonumber\\
&&
\frac{1}{2\pi}\left( \int_{|\omega|\geq 1} d\omega |\widetilde{W}(\omega)|+
\int_{0}^{2\pi} d\phi \int_{0}^{1} dt |\widetilde{W}(t\cos\phi, t\sin\phi)|\right)
\leq
\nonumber\\
&&
\frac{1}{2\pi}(Q+2\pi).\nonumber
\end{eqnarray}
\qed

\begin{lemma}\label{lemma.w2}
For all $W\in\cW(\beta, L)$ we have
\bec
\int_{0}^{\pi} \left| \cR[W](x,\phi)-\cR[W](y,\phi) \right| \, d\phi \leq \frac{Q}{2\pi} |x-y|
\eec
with $Q$ the constant depending on $\beta$ and $L$ defined in Lemma \ref{lemma.w}.
\end{lemma}

\noindent
\textit{Proof.} By (\ref{eq.Fourier transform.tildeW}) and Lemma \ref{lemma.w} we have
\begin{eqnarray*}
&&
\int_{0}^{\pi} \left| \cR[W](x,\phi)-\cR[W](y,\phi) \right| \, d\phi \leq \\
&&
\frac{1}{2\pi}
\int_{0}^{\pi}\int \left|e^{itx}-e^{ity}\right|
\left| \widetilde{W}(t\cos\phi, t\sin\phi) \right| \, d\phi dt \leq
\frac{1}{2\pi} \int_{\R^{2}}\left| \widetilde{W}(w)\right| \,  dw  \leq \frac{Q}{2\pi}|x-y|,
\end{eqnarray*}
where we have used that $\left|e^{itx}-e^{ity}\right| \leq 2|t| |x-y|$.

\qed

\begin{lemma} \label{lem.pbounds}
For all $0<\alpha\leq 1$ and $|x|>1$ there exist constants $c,C$ depending on $\alpha$, such that
\beq
cx^{-(1+2\alpha)}\leq p_{\alpha}(x)\leq Cx^{-(1+2\alpha)}.
\eeq
Moreover there exist constants $C_1$ and $C_2$ such that the first two derivatives of $p_{\alpha}$ satisfy
$ |p_{\alpha}^{\prime}(x)| \leq C_1 x^{-(2+2\alpha)}$, and 
$ |p_{\alpha}^{\prime\prime}(x)| \leq C_2 x^{-(3+2\alpha)}$.

\end{lemma}

\noindent
\textit{Proof} We have
\begin{equation*}
p_\alpha(x,\phi)=\frac{\alpha}{\sqrt{\pi}}\int_0^1 \frac{(1-z)^{\alpha-1/2}}{(1+z)^{1/2}}\exp\left(-x^2\frac{1-z}{1+z}\right)dz,
\end{equation*}
which by the change of variables $u=x\sqrt{\frac{1-z}{1+z}}$ becomes
\begin{equation}\label{eq.palpha}
p_\alpha(x,\phi)=\frac{\alpha 2^{\alpha+1} x}{\sqrt{\pi}}\int_0^x \frac{u^{2\alpha}}{(u^2+x^2)^{\alpha+1}}
\exp(-u^2) \, du.
\end{equation}
It can easily be checked that if $|x|>1$, the right hand side is bounded 
from above  by $Cx^{-(1+2\alpha)}$ and from below by $cx^{-(1+2\alpha)}$, for some positive constants $c,C$.

For the first derivative one can see that
\begin{equation*}
p_\alpha^{\prime}(x) = \frac{1}{x}p_\alpha(x) - \frac{\alpha(\alpha+1) 2^{\alpha+2} x^2}{\sqrt{\pi}} \int_0^x \frac{u^{2\alpha}}{(u^2+x^2)^{\alpha+2}} du + \frac{\alpha}{\sqrt{\pi}x}\exp(-x^2),
\end{equation*}
which can be bounded using the same argument as before to obtain
\begin{equation*}
|p_\alpha^{(1)}(x)| \leq C_1 x^{-(2+2\alpha)},
\end{equation*}
for some $C_1$ and $|x|>1$. For the second derivative the procedure is the same.
\qed

\begin{lemma}\label{lemma.rhoa}
Let $\tau^{a}$ be the diagonal matrix with elements
\beq
\tau^{a}_{k,k}=\frac{1}{4\pi} \int_{0}^{\infty} L_{k} (t^{2}/2) e^{-t^{2}/4}
\frac{t^{2}}
{1+a^{-1}\sinh^{2}\beta t } \, dt.
\eeq
Then
\beq
\tau^{a}_{kk} =O\left( k^{-5/4}(\log a)^{4} \right)
\eeq
\end{lemma}

\noindent
\textit{Proof.} We analyze first the dependance on $a$ for a fixed $k$.
The functions $\{L_{k} (u) e^{-u/2}\}_{k\geq 0}$ form a orthonormal basis of
$L_{2}(\R)$. By using Cauchy-Schwarz followed by Lemma 5 from \cite{Cavalier}
we get
\bec
|\tau^{a}_{k,k}|\leq \frac{1}{4\pi}\left(
\int_{0}^{\infty} \frac{t^{3}}{(1+a^{-1}\sinh^{2}\beta t)^{2} } \, dt \right)^{1/2}\leq
C(\log a)^{3/2},
\eec
for some positive constant $C$.

Let now $a$ be fixed and look at the asymptotic behavior of $\tau^{a}_{k,k}$ as $k\to\infty$. We use the differential equation of the Laguerre polynomials, \cite{Gradshteyn&Ryzhik}  8.979:
\bec
L_{n}(x)=\frac{1}{n}\left( (x-1)L_{n}^{\prime}(x)-xL_{n}^{\prime\prime}(x)\right).
\eec
Thus
\begin{eqnarray}
&&
\frac{d}{dt}L_{n}(t^{2}/2) = t L_{n}^{\prime}(t^{2}/2)
\\
&&
\frac{d^{2}}{dt^{2}}L_{n}(t^{2}/2)= L_{n}^{\prime}(t^{2}/2)+ t^{2}L_{n}^{\prime\prime}(t^{2}/2)
\end{eqnarray}
which implies
\bec
\frac{t^{2}}{2}L_{n}^{\prime\prime}(t^{2}/2)=
\frac{1}{2} \frac{d^{2}}{dt^{2}} L_{n}(t^{2}/2)-\frac{1}{2}t^{-1}\frac{d}{dt}L_{n}(t^{2}/2)
\eec
and
\bec
L_{n}(t^{2}/2)=\frac{1}{2n}\left( (t^{2}-1) t^{-1}\frac{d}{dt}L_{n}(t^{2}/2)- \frac{d^{2}}{dt^{2}} L_{n}(t^{2}/2)\right).
\eec
Using integration by parts we obtain the formula
\bec
\frac{1}{4\pi} \int_{0}^{\infty} L_{k} (t^{2}/2) e^{-t^{2}/4}
\frac{t^{2}}
{1+a^{-1}\sinh^{2}\beta t } \, dt= \frac{1}{k}\int_{0}^{\infty}  L_{k} (t^{2}/2) e^{-t^{2}/4} f(t) \, dt,
\eec
where the function $f$ is given by
\bec
\frac{P_{1}(t)}{1+a^{-1}\sinh^{2}\beta t}+
\frac{a^{-1}(P_{2}(t)\sinh 2\beta t+P_{3}(t)\cosh 2\beta t)}{(1+ a^{-1}\sinh^{2} \beta t)^{2}}+
\frac{P_{3}(t)a^{-2}\sinh^{2} 2\beta t }{(1+a^{-1}\sinh^{2} \beta t)^{3}}.
\eec
with $P_{i}(t)$ polynomials with degree at most four,
whose coefficients do not depend on $a$.

We split the integral into $\int_{0}^{1}$ and $\int_{1}^{\infty}$ and use the following bounds for the behavior of Laguerre polynomials in the two intervals
(see \cite{Szego} Theorem 8.9.12 and Theorem 7.6.4):
\bec
\max_{x\in[1,\infty)} e^{-x/2}|L_{n}(x)| =O( n^{-1/4}),
\eec
and
\bec
L_{n}(x)=x^{-1/4} O(n^{-1/4}),
\eec
uniformly on $(0, 1]$. Thus using Lemma 5 of \cite{Cavalier},
\bec
\left|
\frac{1}{4\pi} \int_{0}^{\infty} L_{k} (t^{2}/2) e^{-t^{2}/4}
\frac{t^{2}}
{1+a^{-1}\sinh^{2}\beta t } \, dt
\right| =O\left( (\log a)^{4}k^{-5/4} \right)
\eec

\qed


%

\begin{lemma}\label{lemma.W0class}
For any $(\beta, L)$ there exists an $\alpha>0$ such that $W_{\alpha}$ belongs to the class $\cW(\beta, L)$.
\end{lemma}

\noindent
\textit{Proof.} By using Minkowski inequality we get
\begin{eqnarray*}
&&
\int e^{2\beta |w|} \left| \widetilde{W}_{\alpha}(w)\right|^{2} \, dw=
2\pi
\int_{0}^{\infty} \left|
\int_{0}^{1}
\sqrt{r}  \frac{f_{\alpha}(z)}{1-z}\exp \left( -r^{2}\frac{1+z}{4(1-z)}+\beta r\right) \, dz
\right|^{2}
\, dr\leq
\\
&&
\left[ \int_{0}^{1} \left(
\int_{0}^{\infty}  r \exp \left( -r^{2}\frac{1+z}{2(1-z)}+2\beta r \right) \, dr
\right)^{1/2}
 \frac{f_{\alpha}(z)}{1-z} \, dz \right]^{2}.
\end{eqnarray*}
The interior integrals satisfies the bound
\begin{eqnarray*}
&&
\int_{0}^{\infty}  r \exp\left( -r^{2}\frac{1+z}{2(1-z)}+2\beta r \right) \, dr \leq C(\beta)(1-z).
\end{eqnarray*}
for some positive constant $C(\beta)$.
Thus
\bec
\int e^{2\beta |w|} \left| \widetilde{W}_{\alpha}(w)\right|^{2} \, dw\leq
C(\beta)\left( \int_{0}^{1} \alpha (1-z)^{\alpha-1/2}\right)^{2}=
C(\beta)\left(\frac{\alpha}{\alpha+1/2}\right)^{2}\to 0,
\eec
as $\alpha\to 0$.

\qed

\begin{lemma} \label{lemma.information}
For $\alpha \leq 1/2$, the Fisher information of the family of densities $\cR W_c^z$ satisfies
\beq
I(c)= 
C^{*} (\log a)^3 R_{0}^{-1}(1+o(1))
\eeq
where $R_{0}$ is defined in \eqref{eq.r0}.
%

\end{lemma}

\noindent
\textit{Proof.} We sketch the proof following the line of \cite{Cavalier} and pointing out where the differences appear. After some transformations, the Fisher information of the family can be brought to the form
\beq \label{eq.infeq}
I(c)= \frac{1}{\pi} \int_0^{\pi} d\phi\int \frac{H_a^2(u)}{\cR [W_c^z]([z,\phi]-u,\phi)} \,du.
\eeq
By expanding $\cR [W_c^z]([z,\phi]-u,\phi)^{-1}$ up to the second order and bounding the second derivative one can show that
\begin{align}
& \left| \int \frac{H_a^2(u)}{\cR [W_c^z]([z,\phi]-u,\phi)} \,du - \frac{1}{\cR [W_c^z]([z,\phi],\phi)} \int H_a^2(u)  \,du \right|  \nonumber \\
& = O\left (\int u^2 H_a^2(u)\, du\right)  = O\left( (\log a)^{2} \right), \quad \text{as} \quad a \to \infty. \label{eq.approx}
\end{align}
Recall that $ \cR [W_{c}^{z}] ([z,\phi], \phi)= \cR [W^{c}] (0, \phi) $ does not depend on 
$\phi$, cf. \eqref{eq.translated}.
Thus we can write
$$
I(c) = \frac{1}{ \cR [W_{c}] (0, \phi)}\int H_{a}^{2} (u)du + O((\log a)^{2})
$$
From  \cite{Cavalier} we have
\begin{align}
& \int H_a^2(u) du = \frac{1}{3\cdot2\cdot(4 \pi \beta)^{2}} (\log a)^{3}(1+o(1))  \quad \text{as} \quad a \to \infty. \label{eq.H_a^2}
\end{align}
From equations (\ref{eq.infeq})--(\ref{eq.H_a^2}) one obtains the desired result.

The main difference in the proof appears when deriving the bound (\ref{eq.approx}). To derive this, one needs a bound of the absolute value of the second derivative 
$$
\left| \frac{\partial^2}{\partial t^2}\frac{1}{\cR [W_c^z](t,\phi)} \right|_{t = [z,\phi]-u}  = 
O(1) , \qquad {\rm as}~~ a\to\infty ,
$$
which is uniform in $u$. We have
\begin{equation} \label{eq.bound_2nd_der}
\left| \frac{\partial^2}{\partial t^2}\frac{1}{\cR [W_c^z](t,\phi)} \right|_{t = [z,\phi]-u} \leq \left|     \frac{p_{c}^{\prime\prime}(u,\phi)}{p_{c}(u,\phi)^2}\right| + \left|\frac{p_c^{\prime}(u,\phi)^2}{p_{c}(u,\phi)^3} \right| 
\end{equation}
where 
$$ 
p_{c}(u,\phi)  = \cR [W_c^z]([z,\phi]-u,\phi) =  \cR [W_c](-u,\phi) =
\cR [W_c](u,\phi)= p_{\alpha} (u,\phi)+ cH_{a}(u),
$$ 
and the derivatives are with respect to the first argument.
If $u$ is in a compact interval, the two terms on the right hand side of 
\eqref{eq.bound_2nd_der} are $O(1)$ as $a\to \infty$ since in that case
$\sup_{|c|<C_{a}} |cH_{a} (u)|\to 0$, and similarly for its derivatives. 
Outside  the interval, we have that $H_{a}$ and its derivatives are exponentially decreasing, and the dominating term is $p_{\alpha}$. Indeed, according to Lemma 
\ref{lem.pbounds}, $p_{\alpha}$  is of order $|u|^{-(1+2\alpha)}$, and its first two derivatives are bounded from above by $|u|^{-(2+2\alpha)}$ and $|u|^{-(3+2\alpha)}$ respectively. 

\qed

\end{document}